\documentclass[12pt,a4paper]{elsarticle}
\usepackage{geometry}
\usepackage{amsmath,amssymb,amsfonts}
\usepackage{mathrsfs}
\usepackage{cases}
\usepackage{chngcntr}
\counterwithout{equation}{section}
\usepackage{color}
\usepackage{ntheorem}

\newtheorem{Th}{Theorem}[section]
\newtheorem{Le}{Lemma}[section]
\newtheorem{De}{Definition}[section]
\newtheorem{Rem}{Remark}[section]
\newtheorem{Cor}{Corollary}[section]
\newtheorem{Ex}{Example}[section]

\newcommand{\bth}{\begin{Th}}
\newcommand{\eeth}{\end{Th}}
\newcommand{\ble}{\begin{Le}}
\newcommand{\eele}{\end{Le}}
\newcommand{\bde}{\begin{De}}
\newcommand{\ede}{\end{De}}
\newcommand{\bre}{\begin{Rem}}
\newcommand{\eere}{\end{Rem}}

\newcommand{\bco}{\begin{Cor}}
\newcommand{\eeco}{\end{Cor}}

\newcommand{\p}{\rm{\bf Proof.}}

\begin{document}
\begin{frontmatter}

\title{Riemann-Hilbert problems for axially symmetric monogenic functions in $\mathbb{R}^{n+1}$}
\author[rvt]{Qian Huang}
\author[rvt]{Fuli He\corref{cor1}}
\author[els1]{Min Ku}
\cortext[cor1]{Corresponding author.\\
 E-mail: huangqianim@163.com, hefuli999@163.com, kumin0844@163.com, 
 }
\address[rvt]{School of Mathematics and Statistics, Central South University, Changsha 410083, China}
\address[els1]{Department of Computing Science, University of Radboud, 6525 EC Nijmegen, Netherlands}

\begin{abstract} We focus on the Clifford-algebra valued variable coefficients Riemann-Hilbert boundary value problems $\big{(}$for short RHBVPs$\big{)}$ for axially monogenic functions on Euclidean space $\mathbb{R}^{n+1},n\in \mathbb{N}$. With the help of Vekua system, we first make one-to-one correspondence between the RHBVPs considered in axial domains and the RHBVPs of generalized analytic function on complex plane. Subsequently, we use it to solve the former problems, by obtaining the solutions and solvable conditions of the latter problems, so that we naturally get solutions to the corresponding Schwarz problems. In addition, we also use the above method to extend the case to RHBVPs for axially null-solutions to $\big{(}\mathcal{D}-\alpha\big{)}\phi=0,\alpha\in\mathbb{R}$.
\end{abstract}

\begin{keyword} Riemann-Hilbert problem, Clifford analysis, axial symmetry, generalized Cauchy-Riemann operator, generalized analytic funcion.
\vspace{0.2cm}

\end{keyword}

\end{frontmatter}

\section{Introduction}

The original version of RHBVPs for analytic functions were was presented by Bernhard Riemann in 1851 \cite{R}. Since Hilbert published his groundbreaking work \cite{H}, in which the problem was associated with singular integral operators, the theory of RHBVPs has started to flourish, see, Refs. \cite{Ga,Lu,Mu,Fo,Ba,Deift,Ch,Ab,Be,WD,KHW}. Besides this, nowadays, the study of such problems have been received wide attention in application fields. This is because the study of many practical problems from mechanics, physics, statistical physics, engineering technology and from theory of orthogonal polynomials and random matrices, see, e.g. \cite{Lu,Mu,Deift,Ch,Ab}, can be attributed to that of a singular integral equation,  see, Refs. \cite{Ga,Mu}, or of a Riemann-Hilbert problem.
Such that, it is of interest to see what it looks like in higher Euclidean dimensions from the view of not only theory values but also practical applications.

Clifford combined the geometric structure and algebraic theory in high-dimensional space and created Clifford algebra, which is an associative but noncommutative algebraic structure. Clifford analysis is the analysis of classical function theory on Clifford algebra. It is the extension of real analysis, complex analysis and quaternion analysis to high-dimensional space. The introduction of Dirac operator in Euclidean space and the appearance of its null-solution monogenic function are important milestones to promote the development of Clifford analysis.
Compared with the classical theory of holomorphic functions in complex analysis, many outstanding scholars have done a lot of research, see, Refs. \cite{BDS,GM,DSS}.
To our knowledge, since 1990s, RHBVPs have been extended into higher dimensional spaces of $\mathbb{R}^{n+1},n\in \mathbb{N}$, see, \cite{ZH}. Afterwards, many scholars are devoting to the higher dimensional RHBVPs in the setting of quaternion analysis or Clifford analysis.
For instance, in \cite{GD,ABP,BD,GZ,KKW,KK,CKK,KFKC} scholars have considered the RHBVPs with constant coefficients and obtained their explicit solutions.
In \cite{HKSB1}, the authors went one step further and solved the variable coefficient RHBVPs for axial monogenic functions in $\mathbb{R}^4$ by using the Fueter theorem in $\mathbb{R}^4$. More details can also be found in \cite{HKSB2,KWHK,KHH}. However, the method there is not applicable to study the variable coefficients RHBVPs for monogenic functions of axial type in $\mathbb{R}^3$, not to mention for axially symmetric monogenic functions defined in higher dimensional Euclidean space $\mathbb{R}^{n+1}$, because the derived boundary value conditions from those of the RHBVPs considered are too complicated. Moreover, besides this, as far as we know, few study of the variable coefficients RHBVPs for monogenic functions, even though for the special cases of axially symmetric monogenic functions defined in higher dimensional Euclidean space, can be seen in references up to date.

Based on all of these facts above, in this paper we mainly investigate Clifford-algebra valued variable coefficients RHBVPs in $\mathbb{R}^{n+1}$ for axially monogenic functions. We first convert them to RHBVPs of the generalized analytic function on the complex plane with the help of Vekua system, after that, we elaborate the form of the general solution of the RHBVPs for generalized analytic function and their solvable conditions. Afterwards, we use these conclusions to solve Clifford-algebra valued variable coefficient RHBVPs, whose boundary date is H\"older continous, for axial monogenic function, in what follows the corresponding Schwarz problems are solved.
Furthermore, we extend the results obtained to null-solutions with axial symmetry to $\big{(}\mathcal{D}-\alpha\big{)}\phi=0,\alpha\in\mathbb{R}$.

The paper proceeds as follows. In section 2, we introduce the definitions and notations of Clifford analysis that we need in what follows. Section 3 is devoted to solving Clifford-algebra valued variable coefficients RHBVPs in $\mathbb{R}^{n+1}$ of the axial monogenic functions. In section 4, we extend the conclusions from the previous section to boundary value problems of Riemann-Hilbert type for null-solutions with axial symmetry to a perturbed generalised Cauchy-Riemann operator.

\section{Preliminaries}

In this section, we in brief review some notions an definitions of Clifford algebra. For more details we refer to \cite{BDS,GM,ZH}.

Let $\big{\{}e_{1}, e_{2},\ldots,e_{n}\big{\}},n\in \mathbb{N}$ be the orthonormal
basis of $\mathbb{R}^n$, satisfying:
$e^{2}_{j}=-1$, if $j=1,2,\ldots,n$,
$e_{i}e_{j}+e_{j}e_{i}=0,$ if $1\leq i<j\leq n$, where $n\in \mathbb{N}$.
For arbitrary element $\alpha$ in real Clifford algebra $\mathbb{R}_{0,n}$ constructed over $\mathbb{R}^n$, it can be expressed as
$$\alpha=\sum_{A\subseteq\mathcal{N}}\alpha_Ae_A:\alpha_A\in\mathbb{R},\mathcal{N}=\{1,2,\ldots,n\}$$
with $e_{A}=e_{r_1}e_{r_2} \ldots e_{r_l}$, $A=\{r_1,r_2,\ldots,r_l\}\subset\mathcal{N}$ satisfying $1\leq r_1<r_2<\ldots<r_l\leq n,$ and when $A=\varnothing,$ we denote $e_{\varnothing}=1.$

The elements of the form $x=x_{0}+\underline{x}=x_0+\sum\nolimits^n_{j=1}x_je_j,$ also called $1$-vectors or parevectors, can be one-to-one mapped to the elements in $\mathbb{R}^{n+1}=\mathbb{R}\oplus \mathbb{R}^n$ through $\sigma:x\mapsto(x_0,x_1,\ldots,x_n).$
Denote the vector and scalar part of $x$ by $\mathrm{Ve}(x)\triangleq\underline{x}$ and $\mathrm{Sc}(x)\triangleq x_0$.

The conjugation $\overline{x}$ of the $1$-vector $x=x_{0}+\underline{x}$ is given as $\bar{x}=x_0+\sum\nolimits^n_{j=1}x_j\bar{e}_j=x_0-\underline{x},$ with $\bar{e}_j=-e_j,j=1,2,\ldots,n,$ so that we have $\overline{xy}=\bar{y}\bar{x}.$ The norm of arbitrary $x\in\mathbb{R}^{n+1}$ is the Euclidean norm $|x|=(\sum\nolimits^{n}_{j=0}\big{|}x_{j}\big{|}^{2})^{\frac{1}{2}}.$ Moreover, if $x\in\mathbb{R}^{n+1}\backslash\{0\},$ then the inverse $x^{-1}$ exists and $x^{-1}:=\bar{x}\cdot|x|^{-2}$, i.e.\ $xx^{-1}=x^{-1}x=1.$ Furthermore, $\mathrm{i}$ is the imaginary unit in $\mathbb{C}$, where $\mathbb{C}$ denotes the set of all complex numbers. $\mathbb{C}_n = \mathbb{R}_{0,n} \otimes\mathbb{C}$. Hence, $ \mathbb{R}_{0,n} \subset \mathbb{C}_n $.

\begin{De}
 Let $\Omega\subset\mathbb{R}^{n+1}$ be a non-empty open subset, for arbitrary $x\in\Omega$, we define the set
\begin{equation}
[x]=\big{\{}y:y=\mathrm{Sc}(x)+\mathcal{I}|\underline{x}|,\mathcal{I}\in S^{n-1}\big{\}},
\end{equation}
 with $S^{n-1}=\big{\{}\underline{x}\subset\mathbb{R}^n:|\underline x|=1\big{\}}$, we say that $\Omega$ is axially symmetric, if the set $[x]$ is contained in $\Omega\subset\mathbb{R}^{n+1}$ for arbitrary $x\in\Omega$.
\end{De}

\begin{Rem}{\rm
Obviously, the  upper half space
$\mathbb{R}^{n+1}_{+}=\big{\{}x\in\mathbb{R}^{n+1}\big{|}x_0>0\big{\}}$ and the ball centered at the origin of $\mathbb{R}^{n+1}$ are both axially symmetric domains.}
\end{Rem}

Let $\Omega$ be an axially symmetric domain of $\mathbb{R}^{n+1}$ whose boundary $\partial\Omega$ is smooth.
All $\mathbb{R}_{0,n}$-valued functions defined on $\Omega$ have the form $\phi=\sum\nolimits_A\phi_{A}e_{A},$ with $\phi_{A}:\Omega\rightarrow\mathbb{R}$. The function $\phi$ is said to be continuous, continuously differentiable, $\mathcal{L}_{p}$-integral, H\"{o}lder continuous if every component $\phi_{A}$ has the property. The corresponding spaces are
denoted by
$\mathcal{C}\big{(}\Omega,\mathbb{R}_{0,n}\big{)}$,
$\mathcal{C}^{1}\big{(}\Omega,\mathbb{R}_{0,n}\big{)}$ ,
$\mathcal{L}_{p}\big{(}\Omega,\mathbb{R}_{0,n}\big{)}\big{(}1<p<+\infty\big{)}$,
$\mathcal{H}^{\mu}\big{(}\Omega,\mathbb{R}_{0,n}\big{)}\big{(}0<\mu\leq1\big{)}$, respectively.

We introduce the concept of generalised Cauchy-Riemann operator
\begin{equation}
\mathcal{D}=\partial_{x_{0}}+\mathcal{D}_{\underline{x}}=\partial_{x_{0}}+\sum\limits^{n}_{j=1}e_{j}\partial_{x_{j}}
\end{equation}
 in Euclidean space $\mathbb{R}^{n+1}$.
It is obvious that $\overline{\mathcal{D}}\mathcal{D}=\sum\nolimits^{n}_{j=0}\partial^{2}_{x_j}=\Delta,$ where $\Delta$ denotes the Laplacian in $\mathbb{R}^{n+1}.$

\begin{De}
 A continuously differentiable $\mathbb{R}_{0,n}$-valued function defined on  $\Omega$ is said to be  (left-)monogenic if and only if $\mathcal {D}\phi=0$. If a function is monogenic and of axial type, it is called axially monogenic. \end{De}

 \begin{Rem}{\rm
 In fact, by Fueter's theorem in Refs. \cite{KFKC,Fu,So,CSS}, axially monogenic functions are functions in the form of
\begin{eqnarray}
\phi(x)=A\big{(}x_{0},r\big{)}+ \underline{\omega} B\big{(}x_{0},r\big{)},\label{axial}
\end{eqnarray}
with $x=x_{0}+\underline{x}=x_{0}+r \underline{\omega} \in\mathbb{R}^{n+1},r=|\underline{x}|,\underline{\omega} \in\{\underline{x}=\sum\nolimits^{n}_{j=1}x_{j}e_{j}:|\underline{x}|=1,x_{1}\ldots,x_{n}\in\mathbb{R}$\}, where $\mathbb{R}$-valued functions $A\big{(}x_{0},r\big{)}$ and $B\big{(}x_{0},r\big{)}$ satisfy a special kind of the Vekua system,
\begin{equation}
\left\{\begin{array}{ll} \partial_{x_0}A-\partial_rB=\frac{n-1}{r}B,\\
\partial_{x_0}B+\partial_rA=0,
\end{array}\right.\label{Vesy}
\end{equation}
with $\partial_{x_0},\partial_r$ denote $\frac{\partial}{\partial_{x_{0}}},\frac{\partial}{\partial_r}$ respectively. Moreover, we define $\mathrm{Re}\ \phi = A,\mathrm{Im}\ \phi = B$.
}
\end{Rem}

In the following context, when the considered functions are defined on $\overline\Omega\subset\mathbb{R}^{n+1}$ with values in $\mathbb{R}_{0,n}$, we mean that they are of axial type if without more explanation.

\begin{Rem}{\rm
Denote by $\mathcal{M}\big{(}\Omega,\mathbb{R}_{0,n}\big{)}$ the set of the monogenic functions of axial type defined in $\Omega$, which actually forms a right-module.}
\end{Rem}

 Moreover, we denote by $D\subset\mathbb{C}_{+}$ the projection of the axially symmetric domain $\Omega\subset\mathbb{R}^{n+1}$ onto the $\big{(}x_{0},r\big{)}$-plane, where $\mathbb{C}_{+}$ is the upper half of the $\big{(}x_{0},r\big{)}$-plane.

\section{Boundary value problems of Riemann-Hilbert type for axially monogenic functions}

In this section we study the Clifford-algebra valued variable coefficients RHBVPs for monogenic functions of axial type defined in Euclidean space $\mathbb{R}^{n+1}$. It is namely to
find an axially monogenic function $\phi\in\mathcal
{C}^{1}\big{(}\Omega,\mathbb{R}_{0,n}\big{)}$, satisfying
\begin{equation}
\left\{\begin{array}{ll} \mathcal {D}\phi(x)=0,&\quad x\in \Omega,\\
\mathrm{Re}\Big{\{}\lambda(t)\phi(t)\Big{\}}=g(t),&\quad t\in
\partial \Omega,
\end{array}\right.\label{prob1}
\end{equation}
where $\lambda(t)=\lambda_1(t)+\underline\omega\lambda_2(t)$ is a $\mathbb{R}_{0,n}$-valued function, with $\lambda_1(t),\lambda_2(t)$ being both $\mathbb{R}$-valued functions defined on $\partial\Omega$, $\mathcal {D}=\partial_{x_{0}}+\sum\nolimits^{n}_{j=1}e_{j}\partial_{x_{j}}$ is the generalised Cauchy-Riemann operator, and $\mathbb{R}$-valued function $g$ is defined on $\partial\Omega$. Hereby, $\lambda,g$ both belong to some function class, for example $\mathcal{H}^{\mu}\big{(}\Omega,\mathbb{R}_{0,n}\big{)}\big{(}0<\mu\leq1\big{)}$ or else.

To do this, we first establish one-to-one correspondence between the considered Riemann-Hilbert problems and Riemann-Hilbert problems for generalised analytic functions on the complex plane, and then we will use it and obtain solutions and the solvable conditions for those who have H\"older continuous boundary value. In addition, we obtain expressions of solutions to Schwarz problems as a special case.

\begin{Rem}{\rm\
When $\lambda_2 \equiv 0$ and $n=3$, the RHBVP \eqref{prob1} reduces to the case of that studied in Refs. \cite{HKSB1,HKSB2,KWHK,KHH}}
\end{Rem}

For this kind of problem, to transfer it to a RHBVP for a generalised analytic function on the complex plane, we will prove the following theorem first.

\begin{Th}\label{th3.1}
The RHBVP \eqref{prob1} is equivalent to finding a holomorphic function $w$ in $D\subset\mathbb{C}_{+}$, satisfying
\begin{equation}
\left\{\begin{array}{ll} \partial_{\bar{z}}w(z)+\frac{(n-1)\mathrm{i}}{4r} w(z)-\frac{(n-1)\mathrm{i}}{4r}\overline{w(z)}=0,&\quad z\in D\subset\mathbb{C}_{+}\\
\mathrm{Re}\Big{\{}\lambda(z)w(z)\Big{\}}=g(z),&\quad z\in
\partial D\subset\mathbb{C}_{+},
\end{array}\right.\label{trans1}
\end{equation}
where
$D\subset\mathbb{C}_{+}$ is the projection of $\Omega\subset\mathbb{R}^{n+1}$, $z=x_{0}+\mathrm{i}r,\partial_{\bar{z}}=\frac{1}{2}\big{(}\partial_{x_{0}}+\mathrm{i}\partial_{r}\big{)},$ and $g:\partial D\rightarrow\mathbb{R}$ is a real-valued function, $\lambda(z)=\lambda_1(z)+\mathrm{i}\lambda_2(z):\partial D\rightarrow\mathbb{C}$ is a complex-valued function.
\end{Th}

\noindent\p \quad We set $\phi(x)=A(x_0,r)+\underline\omega B(x_0,r)$ with $A,B$ as Term \eqref{axial}, and let $w(z)=A(x_0,r)+\mathrm{i} B(x_0,r),z=x_{0}+\mathrm{i}r$. By \eqref{Vesy}, we have
\begin{align*}
\partial_{\bar{z}}w&=\frac12\left(\partial_{x_0}+\mathrm{i}\partial_r\right)(A+\mathrm{i}B)\\
&=\frac12\left(\partial_{x_0}A+\mathrm{i}\partial_{x_0}B+\mathrm{i}\partial_rA-\partial_rB\right)\\
&=\frac{(n-1)B}{2r}=\frac{n-1}{2r}\frac{w-\bar w}{2\mathrm{i}}\\
&=-\frac{(n-1)\mathrm{i}}{4r}(w-\bar w),
\end{align*}
i.e.,
\begin{equation}
\partial_{\bar{z}}w+\frac{(n-1)\mathrm{i}}{4r}w-\frac{(n-1)\mathrm{i}}{4r}\bar w=0.
 \end{equation}
Moreover,
\begin{align*}
\mathrm{Re}\Big{\{}(\lambda_1+\underline\omega\lambda_2)(A+\underline\omega B)\Big{\}}&=\lambda_1 A-\lambda_2 B\\
&=\mathrm{Re}\Big{\{}(\lambda_1+\mathrm{i}\lambda_2)(A+\mathrm{i} B)\Big{\}}.
\end{align*}

This implies that the RHBVP \eqref{prob1} is reduced to
\begin{displaymath}
\left\{\begin{array}{ll} \partial_{\bar{z}}w(z)+\frac{(n-1)\mathrm{i}}{4r} w(z)-\frac{(n-1)\mathrm{i}}{4r}\overline{w(z)}=0,&\quad z\in D\subset\mathbb{C}_{+},\\
\mathrm{Re}\Big{\{}\lambda(z)w(z)\Big{\}}=g(z),&\quad z\in
\partial D\subset\mathbb{C}_{+},
\end{array}\right.
\end{displaymath}
where
$D\subset\mathbb{C}_{+}$ is the projection of $\Omega\subset\mathbb{R}^{n+1}$, and $g:\partial D\rightarrow\mathbb{R}$ is a real-valued function,  $\lambda(z)=\lambda_1(z)+\mathrm{i}\lambda_2(z):\partial D\rightarrow\mathbb{C}$ is a complex-valued function. This finishes the proof. $\hfill\square$

\begin{Rem}{\rm\
Theorem \ref{th3.1} is the key to solve the RHBVPs discussed in this context. Moreover, it also can be applied to solve the RHBVPs appeared in Refs. \cite{HKSB1,HKSB2,KWHK,KHH}.}
\end{Rem}

In fact,  the problem \eqref{trans1} is a RHBVP for generalised analytic function (see, Refs. e.g. \cite{Ve,Ro}), and it is solvable if $\lambda\in \mathcal{H}^{\mu}\big{(}\partial D, \mathbb{C}\big{)},g\in \mathcal{H}^{\mu}\big{(}\partial D,\mathbb{R}\big{)}$. From \cite{Ve,Ro}, we know the solution of problem \eqref{trans1} is representable by
\begin{eqnarray}\label{rep}
w(z)=\Psi(z)\mathrm{e}^{\nu(z)},
\end{eqnarray}
and
\begin{eqnarray}\label{t}
\nu(z)=\frac1\pi\iint\limits_D\left(\frac{(n-1)\mathrm{i}}{4\eta}-\frac{(n-1)\mathrm{i}}{4\eta}\frac{\overline{w(\rho)}}{w(\rho)}\right)
\frac{\mathrm{d}\xi\mathrm{d}\eta}{\rho-z},
\end{eqnarray}
with $\rho=\xi+\mathrm{i}\eta$, the function $\nu\in\mathcal{H}^{\mu}(\overline D,\mathbb{C})$ and function $\Phi\in\mathcal{C}^1(D,\mathbb{C})\cap\mathcal{C}(\partial D,\mathbb{C})$ , which satisfying
\begin{eqnarray}
\mbox{Re}\Big{\{}\lambda_0(z)\Psi(z)\Big{\}}=g(z),\quad z\in
\partial D,\label{rh8}
\end{eqnarray}
with $\lambda_0(z)=\lambda(z)\mathrm{e}^{\nu(z)}\in\mathcal{H}^{\mu}(\partial D,\mathbb{C})$.

For the sake of obtaining the explicit solution to RHBVP \eqref{rh8} for holomorphic function $\Psi(z)$, we first assume that $D=\{z:|z|<1\}$.
Since for any integer $m$ $$\lambda_0(z)z^m=|\lambda_0(z)z^m|\mathrm{e}^{\mathrm{i}(\arg\lambda_0(z)+m\arg z)},$$
so when $z\in\partial D$, we have
\begin{eqnarray}
\lambda_0(z)=\left|\lambda_0(z)\right|z^{-m}\mathrm{e}^{\chi(z)}\mathrm{e}^{-p(z)},\label{eqru}
\end{eqnarray}
where $\chi(z)=p(z)+\mathrm{i}q(z)$ is an analytic function on $D$, whose imaginary part on $\partial D$ is $q(z)=\arg\lambda_0(z)+m\arg z$; $m$ is an integer that makes each branch of $q(z)$ single-valued at $\partial D$.
The function $\chi$ can be expressed by Schwarz integral
\begin{equation}
\chi(z)=\frac{1}{2\pi}\int\limits_{\partial D}q(\tau)\frac{\tau+z}{\tau-z}\frac{\mathrm{d}\tau}{\tau}.
\end{equation}
 Since $q(z)\in\mathcal{H}^\mu(\partial D,\mathbb{C})$, we have $\chi(z)\in\mathcal{H}^\mu(\overline D,\mathbb{C})$.
 Introducing the expression \eqref{eqru} to the boundary condition \eqref{rh8} we obtain
\begin{eqnarray}
\mbox{Re}\Big{\{}z^{-m}\mathrm{e}^{\chi(z)}\Psi(z)\Big{\}}=g_1(z),\quad z\in
\partial D,\label{rh11}
\end{eqnarray}
where $$g_1(z)=\frac{g(z)\mathrm{e}^{p(z)}}{\left|\lambda_0(z)\right|}\in\mathcal{H}^\mu(\partial D,\mathbb{C}).$$

\noindent If $m<0$, then problem \eqref{rh11} implies that
\begin{eqnarray}
\Psi(z)=\frac{z^m\mathrm{e}^{-\chi(z)}}{2\pi\mathrm{i}}\int\limits_{\partial D}g_1(t)\frac{t+z}{t-z}\frac{\mathrm{d}t}{t}+\mathrm{i}c_0z^m\mathrm{e}^{-\chi(z)},
\end{eqnarray}
with $c_0$ is a real constant. Considering the continuity of $\Psi(z)$ we obtain $c_0=0$ and
\begin{eqnarray}
\int_0^{2\pi}g_1\left(e^{\mathrm{i}\theta}\right)\mathrm{e}^{-k\mathrm{i}\theta}\mathrm{d}\theta=0,\quad k=0,\ldots,-m+1.\label{cdt<}
\end{eqnarray}
The above conditions guarantee that $\Psi(z)$ is continuous at $z=0$. Thus, $\Psi(z)$ can be expressed by
\begin{eqnarray}
\Psi(z)=\frac{\mathrm{e}^{-\chi(z)}}{\pi\mathrm{i}}\int\limits_{\partial D}\frac{g_1(t)t^m\mathrm{d}t}{t-z}.\label{slt<}
\end{eqnarray}

\noindent If $m\geq0$, the solution of \eqref{rh11} has the form
\begin{eqnarray}
\Psi(z)=\frac{z^m\mathrm{e}^{-\chi(z)}}{2\pi\mathrm{i}}\int\limits_{\partial D}g_1(t)\frac{t+z}{t-z}\frac{\mathrm{d}t}{t}+\mathrm{e}^{-\chi(z)}\sum_{k=0}^{2m}c_kz^k,\label{slt>}
\end{eqnarray}
where $c_0,c_1,\ldots,c_{2m}$ are complex constants satisfying
\begin{eqnarray}
c_{2m-k}=-\bar c_k,\quad k=0,1,\ldots,m.\label{cdt>}
\end{eqnarray}

Further, if the domain $D$ is a general simply connected domain, there exists a conformal mapping $z=\varphi(\gamma)$ according to the Riemann mapping theorem (see, Refs. e.g. \cite{Con})
\begin{eqnarray}
z=\varphi(\gamma):\ D\longrightarrow \{\gamma:\ |\gamma|<1\},\label{rmmp}
\end{eqnarray}
and its inverse mapping $\gamma=\psi(z)$ is also a conformal mapping
\begin{eqnarray}
\gamma=\psi(z):\ \{\gamma:\ |\gamma|<1\}\longrightarrow D.\label{inrmmp}
\end{eqnarray}
Then the RHBVP \eqref{rh8} on the general simply connected domain $D$ is transformed into the problem on the unit disk $\{\gamma:\ |\gamma|<1\}$
\begin{eqnarray}
\mbox{Re}\Big{\{}\widehat{\lambda}_0(\gamma)\widehat{\Psi}(\gamma)\Big{\}}=\widehat{g}(\gamma),\quad |\gamma|<1,\label{rh19}
\end{eqnarray}
where $\widehat{\Psi}(\gamma)=\Psi[\varphi(\gamma)],\widehat{\lambda}_0(\gamma)=\lambda_0[\varphi(\gamma)],
\widehat{g}(\gamma)=g[\varphi(\gamma)]$, with $\widehat{\Psi}$ is an analytic function in $\{\gamma:\ |\gamma|<1\}$, continuous on its closure, and  $\widehat{\lambda}_0,\widehat{g}$ are also continuous in H\"older sense on $|\gamma|=1$.  From the discussion in the previous paragraph, we can easily get the explicit solution $\widehat{\Psi}(\gamma)$ of RHBVP \eqref{rh19} on the unit disk $\{{\gamma\!:}\ |\gamma|<1\}$, thus, the solution of the RHBVP \eqref{rh8} on general simply connected domain $D$ is
$$\Psi(z)=\widehat{\Psi}[\psi(z)].$$

Based on the above discussion of the RHBVP for generalised analytic function and Theorem \ref{th3.1}, we get the following theorem for original problem \eqref{prob1}.

\begin{Th}\label{th3.2}
Given $\lambda\in \mathcal{H}^{\mu}\big{(}\partial \Omega, \mathbb{R}_{0,n}\big{)},g\in \mathcal{H}^{\mu}\big{(}\partial \Omega,\mathbb{R}\big{)}$, and $D$ is the projection of $\Omega\in\mathbb{R}^{n+1}$ on $\mathbb{C}_{+}$, then the solution to the RHBVP \eqref{prob1} has the form
\begin{eqnarray}
\phi(x)=\mathrm{Re}(w)\big{(}x_{0},|\underline{x}|\big{)}+\underline{\omega}\mathrm{Im}(w)\big{(}x_{0},|\underline{x}|\big{)},\quad x\in\Omega,
\end{eqnarray}
with $\mathrm{Im}(w)$ and $\mathrm{Re}(w)$ denoting the imaginary and real part of $w$, respectively, and $w$ itself is a complex-valued function, given by
\begin{eqnarray} \label{eq23}
w(z)=\Psi(z)\mathrm{e}^{\nu(z)},
\end{eqnarray}
and
\begin{eqnarray} \label{eq24}
\nu(z)=\frac1\pi\iint\limits_D\left(\frac{(n-1)\mathrm{i}}{4\eta}-\frac{(n-1)\mathrm{i}}{4\eta}\frac{\overline{w(\rho)}}{w(\rho)}\right)
\frac{\mathrm{d}\xi\mathrm{d}\eta}{\rho-z},
\end{eqnarray}
where $\rho=\xi+\mathrm{i}\eta$, here
\begin{eqnarray}\label{eq25}
\Psi(z)=\widehat{\Psi}[\psi(z)],
\end{eqnarray}
where $\psi(z)=\gamma$ given by \eqref{inrmmp} is the inverse conformal mapping of $z=\varphi(\gamma)$ given by \eqref{rmmp}.

 If $m\geq0$, $\widehat{\Psi}(\gamma)$ is given by
\begin{eqnarray}
\widehat{\Psi}(\gamma)=\frac{\gamma^m\mathrm{e}^{-\widehat{\chi}(\gamma)}}{2\pi\mathrm{i}}
\int\limits_{|t|=1}\widehat{g}_1(t)\frac{t+\gamma}{t-\gamma}\frac{\mathrm{d}t}{t}+\mathrm{e}^{-\widehat{\chi}(\gamma)}\sum_{k=0}^{2m}c_k\gamma^k,
\end{eqnarray}
where
$${\widehat{g}}_1(\gamma)=\frac{\widehat{g}(\gamma)\mathrm{e}^{\widehat{p}(\gamma)}}{|\widehat{\lambda}(\gamma)\mathrm{e}^{\widehat{\nu}(\gamma)}|},$$
with $\widehat{\lambda}(\gamma)=\lambda[\varphi(\gamma)],\widehat{\nu}(\gamma)=\nu[\varphi(\gamma)],
\widehat{g}(\gamma)=g[\varphi(\gamma)], $ and $m,\widehat{\chi}(\gamma),\widehat{p}(\gamma)$ are given by the relation
\begin{eqnarray}\label{eq100}
\widehat{\lambda}_0(\gamma)=\left|\widehat{\lambda}_0(z)\right|\gamma^{-m}\mathrm{e}^{\widehat{\chi}(\gamma)}\mathrm{e}^{-\widehat{p}(\gamma)},\end{eqnarray}
 on $|\gamma|=1$, where $\widehat{\lambda}_0(\gamma)=\widehat{\lambda}(\gamma)\mathrm{e}^{\widehat{\nu}(\gamma)}$  $\widehat{\chi}(\gamma)=\widehat{p}(\gamma)+\mathrm{i}\widehat{q}(\gamma)$ is an analytic function on $|\gamma|=1$, whose imaginary part on $|\gamma|=1$ is $\widehat{q}(\gamma)=\arg\widehat{\lambda}_0(\gamma)+m\arg\gamma$; $m$ is an integer that makes each branch of $\widehat{q}(\gamma)$ single-valued at $|\gamma|=1$. $c_0,c_1,\ldots,c_{2m}$ are constants satisfying \eqref{cdt>}.

 If $m<0$, $\widehat{\Psi}(\gamma)$ is expressed by
\begin{eqnarray}
\widehat{\Psi}(\gamma)=\frac{\mathrm{e}^{-\widehat{\chi}(\gamma)}}{\pi\mathrm{i}}
\int\limits_{|t|=1}\frac{\widehat{g}_1(t)t^m\mathrm{d}t}{t-\gamma},
\end{eqnarray}
 when and only when
 \begin{eqnarray*}
\int_0^{2\pi}\widehat{g}_1\left(e^{\mathrm{i}\theta}\right)\mathrm{e}^{-k\mathrm{i}\theta}\mathrm{d}\theta=0,\quad k=0,\ldots,-m+1
\end{eqnarray*}
 is fulfilled.
\end{Th}

\noindent\p \quad Applying Theorem \ref{th3.1}, associating with the condition that $D$ is the projection of $\Omega$ on $\mathbb{C}_+$,
the problem \eqref{prob1} is equivalent into the case \eqref{trans1}.
Since $\lambda(x)\in \mathcal{H}^{\mu}\big{(}\partial \Omega, \mathbb{R}_{0,n}\big{)},g(x)\in \mathcal{H}^{\mu}\big{(}\partial \Omega,\mathbb{R}\big{)}$, we have that $\lambda(z)\in \mathcal{H}^{\mu}\big{(}\partial D, \mathbb{C}\big{)},g(z)\in \mathcal{H}^{\mu}\big{(}\partial D,\mathbb{R}\big{)}$. By the statement above, the solution of problem \eqref{trans1} is expressed by the relation
\begin{eqnarray*}
w(z)=\Psi(z)\mathrm{e}^{\nu(z)},\qquad\nu(z)=\frac1\pi\iint\limits_D\left(\frac{(n-1)\mathrm{i}}{4\eta}-\frac{(n-1)\mathrm{i}}{4\eta}\frac{\overline{w(\rho)}}{w(\rho)}\right)
\frac{\mathrm{d}\xi\mathrm{d}\eta}{\rho-z},
\end{eqnarray*}
where $\rho=\xi+\mathrm{i}\eta$, $\Psi(z)=\widehat{\Psi}[\psi(z)]$, where $\psi(z)=\gamma$  given by \eqref{inrmmp} inverse the conformal mapping $z=\varphi(\gamma)$  given by \eqref{rmmp} which mapping $D$ to the unit circular disk $|\gamma|<1$ and $\partial D$ to $|\gamma|=1$.

If $m\geq0$, $\widehat{\Psi}(\gamma)$ is given by
\begin{eqnarray}
\widehat{\Psi}(\gamma)=\frac{\gamma^m\mathrm{e}^{-\widehat{\chi}(\gamma)}}{2\pi\mathrm{i}}
\int\limits_{|t|=1}\widehat{g}_1(t)\frac{t+\gamma}{t-\gamma}\frac{\mathrm{d}t}{t}+\mathrm{e}^{-\widehat{\chi}(\gamma)}\sum_{k=0}^{2m}c_k\gamma^k.
\end{eqnarray}
where ${\widehat{g}}_1(\gamma)=\frac{\widehat{g}(\gamma)\mathrm{e}^{\widehat{p}(\gamma)}}{|\widehat{\lambda}(\gamma)\mathrm{e}^{\widehat{\nu}(\gamma)}|}$,
with $\widehat{\lambda}(\gamma)=\lambda[\varphi(\gamma)],\widehat{\nu}(\gamma)=\nu[\varphi(\gamma)],
\widehat{g}(\gamma)=g[\varphi(\gamma)]$, and $m,\widehat{\chi}(\gamma),\widehat{p}(\gamma)$ are given by the relation
$$\widehat{\lambda}_0(\gamma)=\left|\widehat{\lambda}_0(z)\right|\gamma^{-m}\mathrm{e}^{\widehat{\chi}(\gamma)}\mathrm{e}^{-\widehat{p}(\gamma)}$$
 on $|\gamma|=1$, where $\widehat{\lambda}_0(\gamma)=\widehat{\lambda}(\gamma)\mathrm{e}^{\widehat{\nu}(\gamma)}$,  $\widehat{\chi}(\gamma)=\widehat{p}(\gamma)+\mathrm{i}\widehat{q}(\gamma)$ is an analytic function on $|\gamma|=1$, whose imaginary part on $|\gamma|=1$ is $\widehat{q}(\gamma)=\arg\widehat{\lambda}_0(\gamma)+m\arg\gamma$; $m$ is an integer that makes each branch of $\widehat{q}(\gamma)$ single-valued at $|\gamma|=1$. $c_0,c_1,\ldots,c_{2m}$ are constants satisfying \eqref{cdt>}.

If $m<0$, $\widehat{\Psi}(\gamma)$ is expressed by
\begin{eqnarray}
\widehat{\Psi}(\gamma)=\frac{\mathrm{e}^{-\widehat{\chi}(\gamma)}}{\pi\mathrm{i}}
\int\limits_{|t|=1}\frac{\widehat{g}_1(t)t^m\mathrm{d}t}{t-\gamma},
\end{eqnarray}
 when and only when
\begin{eqnarray*}
\int_0^{2\pi}\widehat{g}_1\left(e^{\mathrm{i}\theta}\right)\mathrm{e}^{-k\mathrm{i}\theta}\mathrm{d}\theta=0,\quad k=0,\ldots,-m+1
\end{eqnarray*}
 is fulfilled.

Therefore, we end up with the solution to the RHBVP \eqref{prob1} has the form
\begin{eqnarray*}
\phi(x)=\mathrm{Re}(w)\big{(}x_{0},|\underline{x}|\big{)}+\underline{\omega}\mathrm{Im}(w)\big{(}x_{0},|\underline{x}|\big{)},\quad x\in\Omega.
\end{eqnarray*}
The result establishes.
$\hfill\square$\\

Let $\lambda\equiv1$, we come to the following conclusion, which is a special case of problem \eqref{prob1}.

\begin{Th}\label{th3.3}
Given $g\in \mathcal{H}^{\mu}\big{(}\partial \Omega,\mathbb{R}\big{)}$, and $D$ is the projection of $\Omega\in\mathbb{R}^{n+1}$ on $\mathbb{C}_{+}$ with boundary $\partial D$, then the solution to the Schwarz problem:
find an axially monogenic function $\phi\in\mathcal
{C}^{1}\big{(}\Omega,\mathbb{R}_{0,n}\big{)}$, satisfying
\begin{equation}
\left\{\begin{array}{ll} \mathcal {D}\phi(x)=0,&\quad x\in \Omega,\\
\mathrm{Re}\big{\{}\phi(t)\big{\}}=g(t),&\quad t\in
\partial \Omega,
\end{array}\right.\label{prob2}
\end{equation}
 has the form
\begin{eqnarray*}
\phi(x)=\mathrm{Re}(w)\big{(}x_{0},|\underline{x}|\big{)}+\underline{\omega}\mathrm{Im}(w)\big{(}x_{0},|\underline{x}|\big{)},\quad x\in\Omega,
\end{eqnarray*}
with $\mathrm{Im}(w)$ and $\mathrm{Re}(w)$ denoting the imaginary and real part of $w$, respectively, and $w$ is a complex-valued function, given by \eqref{eq23} with $\nu(z)$ is represented by \eqref{eq24}. And the expressions $\Psi(z)=\widehat{\Psi}[\psi(z)]$ with $\psi(z)=\gamma$ is given by case as follows.

If $m\geq0$, $\widehat{\Psi}(\gamma)$ is given by
\begin{eqnarray}
\widehat{\Psi}(\gamma)=\frac{\gamma^m\mathrm{e}^{-\widehat{\chi}(\gamma)}}{2\pi\mathrm{i}}
\int\limits_{|t|=1}\widehat{g}_2(t)\frac{t+\gamma}{t-\gamma}\frac{\mathrm{d}t}{t}+\mathrm{e}^{-\widehat{\chi}(\gamma)}\sum_{k=0}^{2m}c_k\gamma^k,
\end{eqnarray}
where $${\widehat{g}}_2(\gamma)=\frac{\widehat{g}(\gamma)\mathrm{e}^{\widehat{p}(\gamma)}}{\left|\mathrm{e}^{\widehat{\nu}(\gamma)}\right|},$$
with $\widehat{\nu}(\gamma)=\nu[\varphi(\gamma)],
\widehat{g}(\gamma)=g[\varphi(\gamma)] $,
 and $m,\widehat{\chi}(\gamma),\widehat{p}(\gamma)$ are given by \eqref{eq100} with $\lambda=1$, $c_0,c_1,\ldots,c_{2m}$ are constants satisfying \eqref{cdt>}.

If $m<0$, $\widehat{\Psi}(\gamma)$ is expressed by
\begin{eqnarray}
\widehat{\Phi}(\gamma)=\frac{\mathrm{e}^{-\widehat{\chi}(\gamma)}}{\pi\mathrm{i}}
\int\limits_{|t|=1}\frac{\widehat{g}_2(t)t^m\mathrm{d}t}{t-\gamma},
\end{eqnarray}
 when and only when
\begin{eqnarray*}
\int_0^{2\pi}\widehat{g}_2\left(e^{\mathrm{i}\theta}\right)\mathrm{e}^{-k\mathrm{i}\theta}\mathrm{d}\theta=0,\quad k=0,\ldots,-m+1
\end{eqnarray*} is fulfilled.
\end{Th}

\begin{Rem}{\rm\
Theorem \ref{th3.2} illustrate a way to solve Clifford-algebra valued variable coefficients RHBVPs for axially monogenic functions in $\mathbb{R}^{n+1}$ like problem \eqref{prob1}.}
\end{Rem}

Problem \eqref{prob1} is the RHBVP for
holomorphic functions on $\mathbb{C}$ when the dimension is 2, for which we refer to \cite{Mu,Fo,Deift}. Moreover, the Schwarz problem for holomorphic functions on $\mathbb{C}$ is the 2-dimensional case of problem \eqref{prob2} when $\lambda\equiv1$ for all $x\in\partial\Omega$.

\section{Boundary value problems of Riemann-Hilbert type for axial meta-monogenic functions}

In this section, we generalise the method used when solving problem \eqref{prob1} to null-solutions to the equation $\big{(}\mathcal{D}-\alpha\big{)}\phi=0,\alpha\in\mathbb{R}$. That is to say, find an axially monogenic function $\phi\in\mathcal
{C}^{1}\big{(}\Omega,\mathbb{R}_{0,n}\big{)}$, satisfying
\begin{equation}
\left\{\begin{array}{ll} \left(\mathcal {D}-\alpha\right)\phi(x)=0,&\quad x\in \Omega,\\
\mbox{Re}\Big{\{}\lambda(t)\phi(t)\Big{\}}=g(t),&\quad t\in
\partial \Omega,
\end{array}\right.\label{prob3}
\end{equation}
where $\alpha$ is understood as $\alpha I$, with $I$ being the identity operator, $\lambda(t)=\lambda_1(t)+\underline\omega\lambda_2(t)$ is a $\mathbb{R}_{0,n}$-valued function on $\partial\Omega$, with $\lambda_1(t),\lambda_2(t)$ are both $\mathbb{R}$-valued functions, $\mathcal {D}=\partial_{x_{0}}+\sum\nolimits^{n}_{j=1}e_{j}\partial_{x_{j}}$ is the generalised Cauchy-Riemann operator, and $\mathbb{R}$-valued function $g$ is defined on $\partial\Omega$.

We start with Theorem \ref{th4.1} given below, which is a promotion of Theorem \ref{th3.2}.

\begin{Th}\label{th4.1}
 Given $\lambda\in \mathcal{H}^{\mu}\big{(}\partial \Omega, \mathbb{R}_{0,n}\big{)},g\in \mathcal{H}^{\mu}\big{(}\partial \Omega,\mathbb{R}\big{)}$, and $D$ is the projection of $\Omega\in\mathbb{R}^{n+1}$ on $\mathbb{C}_{+}$ with boundary $\partial D$, then the solution to the RHBVP \eqref{prob3} has the form
\begin{eqnarray*}
\phi(x)=\mathrm{e}^{\alpha x_0}\Big{(}\mathrm{Re}(w)\big{(}x_{0},|\underline{x}|\big{)}+\underline{\omega}\mathrm{Im}(w)\big{(}x_{0},|\underline{x}|\big{)}\Big{)},\quad x\in\Omega,
\end{eqnarray*}
and $w$ is a complex-valued function, given by \eqref{eq23} with $\nu(z)$ is represented by \eqref{eq24}. And the expressions $\Psi(z)=\widehat{\Psi}[\psi(z)]$ with $\psi(z)=\gamma$ is given by case as follows.

If $m\geq0$, and $\widehat{\Psi}(\gamma)$ is given by
\begin{eqnarray*}
\widehat{\Psi}(\gamma)=\frac{\gamma^m\mathrm{e}^{-\widehat{\chi}(\gamma}}{2\pi\mathrm{i}}
\int\limits_{|t|=1}\widehat{g}_3(t)\frac{t+\gamma}{t-\gamma}\frac{\mathrm{d}t}{t}+\mathrm{e}^{-\widehat{\chi}(\gamma)}\sum_{k=0}^{2m}c_k\gamma^k,
\end{eqnarray*}
where $${\widehat{g}}_3(\gamma)=\frac{\widehat{g}_0(\gamma)\mathrm{e}^{\widehat{p}(\gamma)}}{|\widehat{\lambda}(\gamma)\mathrm{e}^{\widehat{\nu}(\gamma)}|},$$
with $g_0=\mathrm{e}^{-\alpha x_0}g, \widehat{g}_0(\gamma)=g_0[\varphi(\gamma)], \widehat{\lambda}(\gamma)=\lambda[\varphi(\gamma)],\widehat{\nu}(\gamma)=\nu[\varphi(\gamma)]$,
 and $m,\widehat{\chi}(\gamma),\widehat{p}(\gamma)$ are given by \eqref{eq100}, $c_0,c_1,\ldots,c_{2m}$ are constants satisfying \eqref{cdt>}.

If $m<0$, $\widehat{\Psi}(\gamma)$ is expressed by
\begin{eqnarray*}
\widehat{\Psi}(\gamma)=\frac{\mathrm{e}^{-\widehat{\chi}(\gamma)}}{\pi\mathrm{i}}
\int\limits_{|t|=1}\frac{\widehat{g}_3(t)t^m\mathrm{d}t}{t-\gamma},
\end{eqnarray*}
 when and only when
  \begin{eqnarray*}
\int_0^{2\pi}\widehat{g}_3\left(e^{\mathrm{i}\theta}\right)\mathrm{e}^{-k\mathrm{i}\theta}\mathrm{d}\theta=0,\quad k=0,\ldots,-m+1
\end{eqnarray*}is fulfilled.
\end{Th}

\noindent\p \quad Evidently, $(\mathcal{D}-\alpha)\phi=\mathcal{D}(\mathrm{e}^{-\alpha x_0}\phi),\alpha\in\mathbb{R}$, so we have
\begin{eqnarray}
(\mathcal{D}-\alpha)\phi=0{\rm\ \ \Leftrightarrow\ to\ }\mathcal{D}(\mathrm{e}^{-\alpha x_0}\phi)=0,\quad\alpha\in\mathbb{R}.
\end{eqnarray}
Hence, problem \eqref{prob3} reduces to
\begin{eqnarray*}
\left\{\begin{array}{ll} \mathcal {D}(\mathrm{e}^{-\alpha x_0}\phi)=0,&\quad x\in \Omega,\\
\mbox{Re}\Big{\{}\lambda(t)\mathrm{e}^{-\alpha x_0}\phi(t)\Big{\}}=\mathrm{e}^{-\alpha x_0}g(t),&\quad t\in
\partial \Omega,
\end{array}\right.
\end{eqnarray*}
where $g:\partial\Omega\rightarrow\mathbb{R}$  and $\lambda(t)=\lambda_1+\underline\omega\lambda_2:\partial\Omega\rightarrow \mathbb{R}_{0,n}$ is of axial type.

Since $g\in \mathcal{H}^{\mu}\big{(}\partial\Omega,\mathbb{R}\big{)}$, then $g_0\triangleq\mathrm{e}^{-\alpha x_0}g\in \mathcal{H}^{\mu}\big{(}\partial\Omega,\mathbb{R}\big{)}$.

Applying Theorem \ref{th3.2}, with $\lambda$ belongs to $\mathcal{H}^{\mu}\big{(}\partial \Omega, \mathbb{R}_{0,n}\big{)}$, the solution to the RHBVP \eqref{prob3} is represented by
\begin{eqnarray*}
\phi(x)=\mathrm{e}^{\alpha x_0}\Big{(}\mbox{Re}(w)\big{(}x_{0},|\underline{x}|\big{)}+\underline{\omega}\mbox{Im}(w)\big{(}x_{0},|\underline{x}|\big{)}\Big{)},\quad x\in\Omega,
\end{eqnarray*}
with $\mathrm{Im}(w)$ and $\mathrm{Re}(w)$ denoting the imaginary and real part of $w$, respectively, and $w$ is a complex-valued function, given by the relation
\begin{eqnarray*}
w(z)=\Psi(z)\mathrm{e}^{\nu(z)},\qquad\nu(z)=\frac1\pi\iint\limits_D\left(\frac{(n-1)\mathrm{i}}{4\eta}-\frac{(n-1)\mathrm{i}}{4\eta}\frac{\overline{w(\rho)}}{w(\rho)}\right)
\frac{\mathrm{d}\xi\mathrm{d}\eta}{\rho-z},
\end{eqnarray*}
where $\rho=\xi+\mathrm{i}\eta$, $\Psi(z)=\widehat{\Psi}[\psi(z)]$, where $\psi(z)=\gamma$  given by \eqref{inrmmp} inverse the conformal mapping $z=\varphi(\gamma)$  given by \eqref{rmmp} which mapping $D$ to a unit circular disk $|\gamma|<1$ and $\partial D$ to $|\gamma|=1$.

If $m\geq0$, $\widehat{\Psi}(\gamma)$ is given by
\begin{eqnarray*}
\widehat{\Psi}(\gamma)=\frac{\zeta^m\mathrm{e}^{-\widehat{\chi}(\gamma)}}{2\pi\mathrm{i}}
\int\limits_{|t|=1}\widehat{g}_3(t)\frac{t+\gamma}{t-\gamma}\frac{\mathrm{d}t}{t}+\mathrm{e}^{-\widehat{\chi}(\gamma)}\sum_{k=0}^{2m}c_k\gamma^k,
\end{eqnarray*}
where ${\widehat{g}}_3(\gamma)=\frac{\widehat{g}_0(\gamma)\mathrm{e}^{\widehat{p}(\gamma)}}{\left|\widehat{\lambda}(\gamma)\mathrm{e}^{\widehat{\nu}(\gamma)}\right|}$,
with $\widehat{\lambda}(\gamma)=\lambda[\varphi(\gamma)],\widehat{\nu}(\gamma)=\nu[\varphi(\gamma)],
\widehat{g}_0(\gamma)=g_0[\varphi(\gamma)]$,
and $m,\widehat{\chi}(\gamma),\widehat{p}(\gamma)$ are given by \eqref{eq100}, $c_0,c_1,\ldots,c_{2m}$ are constants satisfying \eqref{cdt>}.

If $m<0$, $\widehat{\Psi}(\gamma)$ is expressed by
\begin{eqnarray*}
\widehat{\Psi}(\gamma)=\frac{\mathrm{e}^{-\widehat{\chi}(\gamma)}}{\pi\mathrm{i}}
\int\limits_{|t|=1}\frac{\widehat{g}_3(t)t^m\mathrm{d}t}{t-\gamma},
\end{eqnarray*}
 when and only when
\begin{eqnarray*}
\int_0^{2\pi}\widehat{g}_3\left(e^{\mathrm{i}\theta}\right)\mathrm{e}^{-k\mathrm{i}\theta}\mathrm{d}\theta=0,\quad k=0,\ldots,-m+1
\end{eqnarray*} is fulfilled.
The result follows.$\hfill\square$

\begin{Rem}{\rm\
The two-dimensional case of the problem RHBVP \eqref{prob3} is just the RHBVP for meta-analytic functions on $\mathbb{C}$, which can be found in \cite{Ba,Be,WD,KHW}.}
\end{Rem}

Next we take account of the following problem similar to the problem \eqref{prob2}, that is, finding an axially monogenic function $\phi\in\mathcal
{C}^{1}\big{(}\Omega,\mathbb{R}_{0,n}\big{)}$, satisfying
\begin{equation}
\left\{\begin{array}{ll} \left(\mathcal {D}-\alpha\right)\phi(x)=0,&\quad x\in \Omega,\\
\mbox{Re}\Big{\{}\phi(t)\Big{\}}=g(t),&\quad t\in
\partial \Omega,
\end{array}\right.\label{prob4}
\end{equation}
where $g:\partial\Omega\rightarrow\mathbb{R}$.

\begin{Th}
 Given $g\in \mathcal{H}^{\mu}\big{(}\partial\Omega,\mathbb{R}\big{)}$, then the solution to the RHBVP \eqref{prob4} has the form
\begin{eqnarray*}
\phi(x)=\mathrm{e}^{\alpha x_0}\Big{(}\mathrm{Re}(w)\big{(}x_{0},|\underline{x}|\big{)}+\underline{\omega}\mathrm{Im}(w)\big{(}x_{0},|\underline{x}|\big{)}\Big{)},\quad x\in\Omega,
\end{eqnarray*}
  and $w$ is a complex-valued function, given by \eqref{eq23} with $\nu(z)$ is represented by \eqref{eq24}. And the expressions $\Psi(z)=\widehat{\Psi}[\psi(z)]$ with $\psi(z)=\gamma$ is given by case as follows.

If $m\geq0$,  $\widehat{\Psi}(\gamma)$ is given by
\begin{eqnarray*}
\widehat{\Psi}(\gamma)=\frac{\gamma^m\mathrm{e}^{-\widehat{\chi}(\gamma)}}{2\pi\mathrm{i}}
\int\limits_{|t|=1}\tilde{g}_4(t)\frac{t+\gamma}{t-\gamma}\frac{\mathrm{d}t}{t}+\mathrm{e}^{-\widehat{\chi}(\gamma)}\sum_{k=0}^{2m}c_k\gamma^k,
\end{eqnarray*}
where
$${\widehat{g}}_4(\zeta)=\frac{\widehat{g}_0(\gamma)\mathrm{e}^{\widehat{p}(\gamma)}}{|\mathrm{e}^{\widehat{\nu}(\gamma)}|},$$
with $g_0=\mathrm{e}^{-\alpha x_0}g, \widehat{g}_0(\gamma)=g_0[\varphi(\gamma)], \widehat{\lambda}(\gamma)=\lambda[\varphi(\gamma)],\widehat{\nu}(\gamma)=\nu[\varphi(\gamma)]$,
and $m,\widehat{\chi}(\gamma),\widehat{p}(\gamma)$ are given by \eqref{eq100} with $\lambda=1$, $c_0,c_1,\ldots,c_{2m}$ are constants satisfying \eqref{cdt>}.

If $m<0$, $\widehat{\Psi}(\gamma$ is expressed by
\begin{eqnarray*}
\widehat{\Phi}(\gamma)=\frac{\mathrm{e}^{-\widehat{\chi}(\gamma)}}{\pi\mathrm{i}}
\int\limits_{|t|=1}\frac{\widehat{g}_4(t)t^m\mathrm{d}t}{t-\gamma},
\end{eqnarray*}
 when and only when
 \begin{eqnarray*}
\int_0^{2\pi}\widehat{g}_4\left(e^{\mathrm{i}\theta}\right)\mathrm{e}^{-k\mathrm{i}\theta}\mathrm{d}\theta=0,\quad k=0,\ldots,-m+1
\end{eqnarray*} is fulfilled.
\end{Th}

\begin{Rem}{\rm\
With the H\"older continuous boundary data, problems \eqref{prob1}, \eqref{prob2}, \eqref{prob3} and \eqref{prob4} have discussed the Clifford-algebra valued variable coefficients RHBVPs for axial monogenic functions in Euclidean space $\mathbb{R}^{n+1}$. After further observation, we can point out that Problem \eqref{prob1} is the special case of Problem \eqref{prob3} when $\alpha=0$, while problem \eqref{prob4} reduces to Problem \eqref{prob2} when $\alpha$ equals to 0. In addition, for the RHBVPs of axially monogenic functions in $\mathbb{R}^{n+1}$, if we replace the H\"older continuous boundary value condition with $\mathcal{L}_{p}$-integral boundary value condition, with $1<p<+\infty$, all results above can be established following the same approch. Moreover, the Clifford-algebra valued variable coefficients RHBVPs for monogenic functions in Euclidean space $\mathbb{R}^{n+1}$ is still in progress, and will be discussed in the forthcoming paper.}
\end{Rem}

\section*{Acknowledgements}
We would like to thank Prof. Uwe Kaehler for sharing his idea without reservation when we started to study RHBVPs with variable
coefficients for monogenic functions in high-dimensional Euclidean space several years ago.


\section*{Data Availability Statements}
Data sharing not applicable to this article as no datasets were generated or analysed during the current study.

\section*{Conflict of interest}
The authors declare that they have no conflict of interest.



\begin{thebibliography}{1}
\bibitem{R}  Riemann B. Grundlagen f\"{u}r eine allgemeine Theorie der Functionen einer ver\"{a}nderlichen complexen Gr\"{o}sse. G\"{o}ttingen: Adalbert Rente, 1851.
\bibitem{H}  Hilbert D. \"{U}ber eine Anwendung der Integralgleichungen auf ein Problem Funktionentheorie, Verhandl. der III Int. Math. Kongr., 1904, 223-240.
\bibitem{Ga}  Gakhov FD. Boundary Value Problems. Oxford: Pergamon, 1966.
\bibitem{Lu}  Lu JK. Boundary value problems for analytic functions. Singapore: World Scientific, 1993.
\bibitem{Mu} Muskhelishvili NI. Singular integral equations. Leyden: Noordhoff, 1977.
\bibitem{Fo} Fokas AS. A unified approach to boundary value problems. Cambridge: University of Cambridge, 2008.
\bibitem{Ba} Balk MB. On poly-analytic functions. Berlin: Akademie Verlag, 1991.
\bibitem{Deift} Deift P. Orthogonal polynomials and random matrices: a Riemann-Hilbert approach. Providence, Rhode Island:  American Mathematical Society, 2000.
\bibitem{Ch} Chelkak D, Smirnov S. Universality in the 2D Ising model and conformal invariance of fermionic observables. Invent. Math. 2012, 189(3): 515-580.
\bibitem{Ab}  Abreu LD, Feichtinger HG. Function spaces of poly-analytic functions, Harmonic and Complex Analysis and its Applications, 2014(Trends in Mathematics), 1-38.
\bibitem{Be} Begehr H, Schmersau D. The Schwarz problem for poly-analytic functions. ZAA. 2005, 24(2): 341-351.
\bibitem{BDS} Brackx F, Delanghe R, Sommen F. Clifford analysis res notes math. Vol. 76, London: Pitman, 1982.
\bibitem{WD} Wang Y, Du J. Mixed boundary value problems with a shift for a pair of meta-analytic and analytic functions. J. Math. Anal. Appl., 2010, 369(2): 510-524.
\bibitem{KHW}
Ku M, He FL, Wang Y. Riemann Hilbert problems for Hardy space of meta-analytic functions on the unit disc. Complex Anal. Oper. Th., 2018, 12: 457-474.

\bibitem{DSS} Delanghe R, Sommen F, Souc\'{e}k V. Clifford algebra and spinor-valued functions. Dordrecht: Kluwer Academic, 1992.
\bibitem{GD} Gong YF, Du JY. A kind of Riemann and Hilbert boundary value problem for left monogenic functions in
$\mathbb{R}^{m}(m\geq2)$. Complex Variables Theory Appl., 2004, 49(5): 303-318.
\bibitem{ABP} Abreu Blaya R, Bory Reyes J, Pe\~na-Pe\~na D. Jump problem and removable singularities for monogenic functions. J. Geom. Anal. Theory Appl. 2007, 17(1): 1-13.
\bibitem{BD} Bu YD, Du JY. The RH boundary value problem for the $k$-monogenic functions.  J. Math. Anal. Appl., 2008, 347: 633-644.
\bibitem{GZ} G\"urlebeck K, Zhang ZX. Some Riemann boundary value problems in Clifford analysis. Math. Methods Appl. Sci., 2010, 33: 287-302.
\bibitem{KKW} Ku M, K\"ahker U, Wang DS. Riemann boundary value problems on the sphere in Clifford analysis. Adv. Appl. Clifford Algebra, 2012, 22(2): 365-390.
\bibitem{KK} Ku M, K\"ahker U. Riemann boundary value problems on half space in Clifford analysis. Math. Methods Appl. Sci., 2012, 35(18): 2141-2156.
\bibitem{CKK} Cerejeiras P, K\"ahker U, Ku M. On the Riemann boundary value problem for null solutions to iterated generalized Cauchy-Riemann operator in Clifford analysis. Results Math. 2013, 63(3-4): 1375-1394.
\bibitem{KFKC} Ku M, Fu YX, K\"ahker U, Cerejeiras P. Riemann boundary value problems for iterated Dirac operator on the ball in Clifford analysis. Complex Anal. Oper. Th., 2013, 7(3): 673-693.
\bibitem{Fu} Fueter R. Die Funktionentheorie der Differentialgleichungen $\Delta u$=0 und $\Delta\Delta u$=0 mit vier reellen Variablen. Comm. Math. Helv., 1934, 7: 307-330.
\bibitem{So} Sommen F. On a generalization of Fueter's theorem. ZAA., 2000, 19: 899-902.

\bibitem{GM} Gilbert J, Murray M. Clifford algebras and Dirac operators in harmonic analysis. Cambridge: University of Cambridge, 1991.
\bibitem{CSS} Colombo F, Sabadinia I, Sommen F. The Fueter mapping theorem in integral form and the $\mathcal{F}$-functional calculus. Math. Methods Appl. Sci., 2010, 33: 2050-2066.

\bibitem{HKSB1} He FL, Ku M, K\"ahker U, Sommen F, Bernstein S. Riemann-Hilbert problems for monogenic functions in axially symmetric domains. Bound. Value Probl., 2016, 22: 1-11.

\bibitem{HKSB2} He FL, Ku M, K\"ahker U, Sommen F, Bernstein S. Riemann-Hilbert problems for null-solutions to iterated generalized Cauchy-Riemann equations in axially symmetric domains.
    Comput. Math. Appl., 2016, 71(10): 1990-2000.
\bibitem{ZH} Xu ZY, Zhou ZP. On boundary value problems of Riemann-Hilbert type for monogenic functions in a half space of $\mathbb{R}^{m}\left(m\geq2\right)$. Complex Var. Elliptic Equ., 1993, 22(3-4): 181-193.
\bibitem{KWHK} Ku M, Wang Y, He FL, K\"ahker U. Riemann-Hilbert Problems for Monogenic Functions on Upper Half Ball of $\mathbb{R}^4$ . Adv. Appl. Clifford Algebra, 2017, 27(3): 2493-2508.

\bibitem{KHH} Ku M, He FL, He XL. Riemann-Hilbert problems for null-solutions to iterated generalized Cauchy-Riemann equation on upper half ball. Complex Var. Elliptic Equ., 2019, 65(11): 1902-1918.

\bibitem{Ve}  Vekua IN. Generalized analytic functions. London: Pergamon, 1962.

\bibitem{Ro}  Rodin YL. Generalized analytic functions on Riemann surfaces. Berlin: Springer, 1987.

\bibitem{Con}  Conway JB. Functions of one complex variable. New York: Springer-Verlag, 1978.






\end{thebibliography}
 \end{document}